\documentclass[reqno,a4paper,final]{amsart}
\usepackage{amssymb,amstext,amsthm,eucal}
\usepackage{bbold}
\usepackage{array}
\usepackage[latin1]{inputenc}
\renewcommand{\d}{\partial}

\newcommand{\SO}{\mathrm{SO}}

\newcommand{\Spin}{\mathrm{Spin}}
\newcommand{\SP}{\mathrm{Sp}}

\newcommand{\SU}{\mathrm{SU}}

\newcommand{\RR}{\mathbb{R}}
\newcommand{\CC}{\mathbb{C}}
\newcommand{\HH}{\mathbb{H}}

\newtheorem{THEO}{\bf Theorem}

\newtheorem{PR}{\bf Proposition}

\newtheorem{CO}{\bf Corollary}
\newtheorem{LM}{\bf Lemma}

\newcommand{\MUNCH}[1]{\relax}
\allowdisplaybreaks[1]
\begin{document}
\begin{sloppypar}
\title{Twistor spinors with zero on Lorentzian 5-space}
\author{Felipe Leitner}
\address{Institut f{\"u}r Geometrie und Topologie, Universit{\"a}t 
Stuttgart, Pfaffenwaldring 57, D-70569, Germany}
\email{leitner@mathematik.uni-stuttgart.de}
\thanks{The research to this work was supported by a 
Junior Research Fellowship grant of the Erwin-Schr{\"o}dinger-International (ESI) Institute in Vienna 
financed by the Austrian 
Ministry
of Education, Science and Culture (BMBWK) and by a fellowship grant of the 
Sonderforschungsbereich 647 'Space-Time-Matter - Geometric and Analytic Structures' 
of the Deutsche Forschungsgemeinschaft (DFG) at Humboldt University Berlin}
\date{February 2006}

\begin{abstract} 
We present in this paper a $C^1$-metric on an open neighbourhood of the origin
in $\RR^{5}$. The metric is of Lorentzian signature $(1,4)$ and admits a solution 
to the twistor equation for spinors with a unique isolated zero at the origin.
The metric is not conformally flat in any neighbourhood of the origin.
The construction is based on the Eguchi-Hanson metric with parallel 
spinors on Riemannian $4$-space.\\

MSC-class: 53C50, 53A30, 53C27, 58J70, 57S20.
\end{abstract}

\maketitle
\tableofcontents

\section{Introduction}
\label{ab1}
For spinor fields of suitable weight on semi-Riemannian manifolds there exist two conformally covariant linear differential equations
of first order, the Dirac equation and the twistor equation. The twistor equation is overdetermined and the existence of solutions,
which are called conformal Killing spinors
or simply twistor spinors, is constraint by curvature conditions on the underlying space. 
The twistor equation was first introduced by R. Penrose in General
Relativity. In the second half of the 1980th A. Lichnerowicz started a systematic investigation of twistor spinors on 
Riemannian spin manifolds in the context of conformal differential geometry (cf. \cite{Lic88}, \cite{Lic89}, \cite{Lic90}). 
Special solutions of the twistor equation are Killing spinors and parallel spinors, for which
nowadays many geometric structure results are known.

From the view point of conformal geometry twistor spinors with zeros are of particular interest (cf. \cite{Lic90}, \cite{KR95}, 
\cite{KR96}, \cite{KR98}). 
This is for various reasons. In the Riemannian case the length square of a twistor gives rise to a
rescaled Ricci-flat metric in the conformal class on the complement of the zero set, which consists of isolated
points. Such spaces are sometimes called almost conformally Einstein manifolds (cf. \cite{Gov04}). A result 
by A. Lichnerowicz states that a compact Riemannian space admitting a twistor spinor with zero is isometric
to the standard $n$-sphere $S^n$, which is a conformally flat space. Any twistor on the $n$-sphere admits exactly one 
isolated zero. However, a construction by W. K{\"u}hnel and H.-B. Rademacher shows that there exist twistors with 
zeros on complete non-compact Riemannian spaces, which are not conformally flat. Such solutions occur typically on the conformal
completion space to infinity of asymptotically Euclidean spaces with special holonomy (cf. \cite{KR96}, \cite{KR98}). 

In the Lorentzian setting solutions of the twistor equation always give rise to non-trivial conformal Killing vector fields and 
the zero sets of a twistor and its corresponding conformal Killing vector coincide. There are two types of twistor spinors with
zeros. In the first case the associated conformal Killing vector is almost everywhere timelike and the zero set consists of 
isolated points. Outside of the lightcones of the zeros the geometry is locally conformally equivalent to a static monopole based
on a Riemannian space with parallel spinor. In the other case the associated field is lightlike (or zero), 
the zero set consists of disjoint lightlike geodesics and on the complement of the zero set the geometry is locally conformally
equivalent to a generalised pp-wave (cf. \cite{Lei99}, \cite{Lei01}, \cite{Lei04}). 
  
The associated conformal Killing vector to a twistor with zero has the interesting property that 
its (local) flow consists of essential conformal transformations, i.e., transformations which
are not isometries for any metric in the conformal class of the underlying space. 
Essential conformal transformation groups are non-compact. A conjecture by A. Lichnerowicz states that
compact Lorentzian spaces with essential conformal transformation group are conformally flat. 
In the Riemannian case it is well known that the only complete spaces with essential (i.e. non-compact) conformal
transformation group are the $n$-sphere $S^n$ (compact case) and the Euclidean space (non-compact case)
(cf. \cite{Ale72}, \cite{Yos75}, \cite{Oba71}, \cite{LF71}).   
 
In this paper we construct a family of Lorentzian metrics on (non-compact) open subsets of $\RR^5$, which  admit twistor spinors with a unique 
isolated zero. The conformal geometry of these metrics is not flat in any neighbourhood of the zero. The construction is 
based on the conformal completion of the Eguchi-Hanson metric (cf. \cite{KR96}). 
The constructed family of metrics is of class $C^1$, i.e., with respect to the standard coordinates
on $\RR^5$ the coefficients of the metrics are continuously differentiable exactly once. 
The course of the paper is very simple. In section 
\ref{ab2} we present the family
of metrics in question with some chosen frame and a spinor. In section \ref{ab3} we calculate 
that the metrics in the family are of class $C^1$   
and that the given spinor solves the twistor equation.

\section{Metric with frame and spinor}
\label{ab2}
Let us consider the 5-dimensional real vector space $\RR^5$ with canonical coordinates
$x=(x_0,x_1,x_2,x_3,x_4)$. We set $n=5$. The Minkowski metric is given by
\[g_0=-dx_0^2+dx_1^2+dx_2^2+dx_3^2+dx_4^2\ .\]
This metric is of Lorentzian signature $(1,4)$ and is flat on $\RR^5$. We aim to rewrite the Minkowski metric
in cylindrical coordinates. 
So let $E$ be the 4-dimensional
vector subspace in $\RR^5$ defined by $x_0=0$ and 
denote by  \[r=\sqrt{x_1^2+x_2^2+x_3^2+x_4^2}\] the radial coordinate on $E$.
The space $E\smallsetminus \{0\}$ (with deleted origin) is diffeomorphic to $\RR_+\times S^3$. Thereby, 
$S^3$ denotes the $3$-sphere, which is given in $E$ by the equation $r=1$. 
As the group of elements with unit length in $E\cong \HH$ the $3$-sphere $S^3$ is isomorphic to the 
semisimple Lie group $\SU(2)$. 
The round metric $g_{S^3}$ on $S^3$ is $\SU(2)$-invariant and there exist left-invariant $1$-forms
$\sigma_1,\sigma_2$ and $\sigma_3$ on $\SU(2)$ such that 
\[
g_{S^3}\ =\ \sigma_1^2+\sigma_2^2+\sigma_3^2\ .
\]
On $S^3$ in $E$ these left-invariant forms are explicitly given by
\[
\begin{array}{lcl}
\sigma_1&=& \frac{1}{r^2}(- x_2dx_1+x_1dx_2-x_4dx_3+x_3dx_4),\\[2mm]
\sigma_2&=& \frac{1}{r^2}(- x_3dx_1+x_4dx_2+x_1dx_3-x_2dx_4),\\[2mm]
\sigma_3&=& \frac{1}{r^2}(- x_4dx_1-x_3dx_2+x_2dx_3+x_1dx_4)\ .
\end{array}
\]
We denote the dual orthonormal frame on $TS^3$ by 
$\{\frac{\d}{\d \sigma_1},\frac{\d}{\d \sigma_2},\frac{\d}{\d \sigma_3} \}$.
Eventually, we see that the Minkowski metric on $\RR^5\smallsetminus \{r=0\}$ is given 
in cylindrical coordinates by
\[
g_0\ =\ -dx_0^2+dr^2+r^2(\sigma_1^2+\sigma_2^2+\sigma_3^2)\ .
\] 
We know that this metric can be smoothly completed to the singular set $\{r=0\}$
of the cylindrical coordinate system (which is a real line in $\RR^5$).
The result is the   
Minkowski metric $g_0$ on $\RR^5$.

Now let us define the cone
\[
L\ :=\ \{\ (x_0,x_1,x_2,x_3,x_4)\in\RR^5:\ r\leq|x_0|\ \} 
\]
with singular point at the origin of $\RR^5$. The boundary set of the cone $L$ in $\RR^5$ is
\[
L_o\ :=\ \{\ (x_0,x_1,x_2,x_3,x_4)\in\RR^5:\ r=|x_0|\ \}\ .
\]
As next we define the radial coordinate
\[ r_o\ :=\ \left\{\begin{array}{ccl}
0&&\ \mbox{on}\ \ L\\[2mm]
\frac{r^2-x_0^2}{r}&&\ \mbox{on}\ \
\RR^5\smallsetminus L
\end{array}\right.\quad .
\]
Furthermore, let $a>0$ be a real parameter.
Then we set 
\[
B_a\ :=\ \{\ (x_0,x_1,x_2,x_3,x_4)\in\RR^5:\ 0<r_o<\frac{1}{a}\ \}
\]
and $\tilde{B}_a:= B_a\cup L$. Both sets $B_a$ and $\tilde{B}_a$
are open in $\RR^5$ for all $a>0$. The set $B_a$ is a subset of $\RR^5\smallsetminus L$ and $\tilde{B}_a$ is simply connected.
We also denote the set $B_a^>:=\tilde{B}_a\smallsetminus \{r=0\}$, where the real line $\{r=0\}$ is deleted.   

On $\tilde{B}_a$ we define a family of pointwise symmetric bilinear forms $g_a$, $a>0$, as follows. Let
\[ g_a\ :=\ \left\{\begin{array}{lcl}\ \
g_0\ -\ r^2(ar_o)^4\cdot\sigma_3^2\ +\ a^4(r\beta)^{-2}r_o^2\cdot \alpha^2
&&\ \mbox{on}\ \ \ B^>_a\\[2mm]
\ \ g_0&&\ \mbox{on}\ \ \{r=0\}\\[2mm]
\end{array}\right.\quad ,
\]
whereby we set 
\[ \beta:=\sqrt{1-(ar_o)^4}\qquad\mbox{and}\qquad \alpha:= (r^2+x_0^2)dr\ - \
2x_0r dx_0\ .\] 
Obviously, the symmetric bilinear form $g_a$ is smoothly defined on $\tilde{B}_a\smallsetminus L_o$ for all $a>0$, and
by definition, $g_a$ restricted to $L\smallsetminus L_o$ is the flat Minkowski metric.
The symmetric bilinear form $g_a$ can be rewritten on $B^>_a$ as
\[g_a\ =\ -dx_0^2\ +\ dr^2\ +\ r^2(\ \sigma_1^2+\sigma_2^2+\beta^2\sigma_3^2\ )\ +\ 
a^4(r\beta)^{-2}r_o^2\cdot \alpha^2\quad .\]

\begin{PR}\label{PR1}
The symmetric bilinear form $g_a$ is a $C^1$-metric of Lorentzian signature $(1,4)$ on the subset $\tilde{B}_a$ of $\RR^5$
for all $a>0$.
The metric $g_a$ is not of class $C^2$. 
\end{PR}

We want to do some geometric discussion of the Lorentzian metric $g_a$. 
The restriction of $g_a$ to the disk $E\cap B_a$ (with deleted origin) in $E\cong\RR^4$ is given by
\[
h_a\ :=\ \frac{dr^2}{1-(ar)^4}\ +\ r^2(\ \sigma_1^2+\sigma_2^2+(1-(ar)^4)\sigma_3^2\ )\ .
\]
This is a Riemannian metric on $E\cap B_a$, which admits a smooth (even analytic) extension
to the origin. Off the origin, the metric $g_a$ is conformally equivalent to the
Eguchi-Hanson metric 
\[g_{EH}\ :=\ \frac{dR^2}{1-(a/R)^4}\ +\ R^2\left(\ \sigma_1^2+\sigma_2^2+(1-(a/R)^4)\sigma_3^2\ \right)\ 
.\]
In fact, with $R:=1/r$\ on\ $E\smallsetminus \{0\}$ it holds $g_{EH}=\frac{1}{r^4}\cdot h_a$. 
The Eguchi-Hanson metric is an asymptotically Euclidean hyperk{\"a}hler metric with irreducible 
holonomy group $\SU(2)=\SP(1)$ and admits a $2$-dimensional space of parallel spinors.
In particular, $g_{EH}$ is Ricci-flat, but it is not conformally flat. In fact, $g_{EH}$ is half-conformally
flat, i.e., the Weyl curvature tensor $W=W^++ W^-$ is anti-selfdual ($W^+=0$) (cf. \cite{EH78}, \cite{KR96}). 
We set
\[
\tilde{g}_a:=\frac{1}{(r^2-x_0^2)^2}\cdot g_a\quad\qquad \mbox{on}\ \ \tilde{B}_a\smallsetminus L_o\ .
\]

\begin{PR}\label{PR2} a) Let $\tilde{g}_a=(r^2-x_0^2)^{-2}\cdot g_a$, $a>0$, be a conformally equivalent metric 
to $g_a$ on $\tilde{B}_a\smallsetminus L_o$. 
Then 
\begin{enumerate}
\item the metric $\tilde{g}_a$ is flat for $r<|x_0|$.
\item For $|x_0|<r$ it holds 
\[\tilde{g}_a\ =\ -ds^2\ +\ g_{EH}\ ,\] 
whereby\ $s:=\frac{-x_0}{r^2-x_0^2}$\ and\ $R:=\frac{r}{r^2-x_0^2}$\ .
\end{enumerate}
In particular, $\tilde{g}_a$ is Ricci-flat on $\tilde{B}_a\smallsetminus L_o$.

b) The Weyl tensor $W^{g_a}$ of the smooth Lorentzian metric $g_a$ on $\tilde{B}_a\smallsetminus L_o$ admits
a continuous extension of class $C^1$ to the singular set $L_o$. For this extension it holds $W^{g_a}\equiv 0$ on $L$ and  $W^{g_a}\neq 0$ on $B_a$, 
i.e.,  
$g_a$ is not conformally flat.
\end{PR}
We note that the Ricci-curvature tensor of the metric $g_a$ on $\tilde{B}_a\smallsetminus L_o$ does not admit a continuous 
extension to $\tilde{B}_a$. 
With $\mu:=\ln|r^2-x_0^2|$ it is 
\[
Ric^{g_a}\ =\ -(n-2)\cdot(\ Hess^{\tilde{g}_a}(\mu)\ -\ d\mu^2 \ )\ -\ (\ \Delta^{\tilde{g}_a}\mu\ +\ (n-2)\cdot|d\mu|^2 \ )\cdot 
\tilde{g}_a\ .
\]
Furthermore, we note that the hypersurface $\{s=0\}$ is totally geodesic with 
respect to the metric\ $-ds^2+g_{EH}$.
This implies that the disk $E\cap \tilde{B}_a$  
is a totally umbilic hypersurface in
$(\tilde{B}_a,g_a)$.  

As next we define on $B_a^>$ with metric 
$g_a$ 
an orthonormal 
frame $e=\{e_0,e_1,e_2,e_3,e_4\}$ 
in the following way. Let 
\[
T\ :=\ -(r^2+x_0^2)\frac{\d}{\d r}\ -\ 2rx_0\frac{\d}{\d x_0}
\]
be a vector field on $\RR^5$. We set 
\[\begin{array}{lcl}
e_0&:=&
\frac{\partial}{\partial x_0}\ -\ \frac{2x_0}{r}\cdot\frac{a^4r_o^2}{1+\beta}\cdot T\\[6mm]
e_1&:=&
\frac{\partial}{\partial r}\ +\ \frac{r^2+x_0^2}{r^2}\cdot\frac{a^4r_o^2}{1+\beta}\cdot T\\[6mm]
e_2&:=&r^{-1}\cdot\frac{\d}{\d \sigma_1}\\[6mm]
e_3&:=&r^{-1}\cdot\frac{\d}{\d \sigma_2}\\[6mm]
e_4&:=&(r\beta)^{-1}\cdot\frac{\d}{\d \sigma_3}\quad .
\end{array}\]
\begin{LM}\label{LMH}
The orthonormal frame $e=\{e_0,e_1,e_2,e_3,e_4\}$ on $B^>_a$ 
is of class $C^1$. 
\end{LM}

We proceed by introducing spinor calculus on $\tilde{B}_a$ (cf. e.g. \cite{Baum81}). Let $\Spin(1,4)$ denote the spin group 
with universal 
covering map
$\lambda: \Spin(1,4)\to SO(1,4)$ onto the special orthonormal group and let $Cl(1,4)$ be the Clifford algebra. The 
complex spinor module $\Delta_{1,4}$ is 
isomorphic to $\CC^4$ and a realisation
of the Clifford algebra $Cl(1,4)$ on $\Delta_{1,4}\cong\CC^4$ is given by
\[
\gamma_0=
\left(
\begin{array}{rrrr}
-1&0&0&0\\
0&-1&0&0\\
0&0&1&0\\
0&0&0&1
\end{array}
\right),\qquad
\gamma_1=
\left(
\begin{array}{rrrr}
0&0&-1&0\\
0&0&0&1\\
1&0&0&0\\
0&-1&0&0
\end{array}
\right),\]
\[
\gamma_2=
\left(
\begin{array}{rrrr}
0&0&-i&0\\
0&0&0&-i\\
-i&0&0&0\\
0&-i&0&0
\end{array}
\right),\qquad
\gamma_3=
\left(
\begin{array}{rrrr}
0&0&0&-1\\
0&0&-1&0\\
0&1&0&0\\
1&0&0&0
\end{array}
\right),\]
\[
\gamma_4=
\left(
\begin{array}{rrrr}
0&0&0&-i\\
0&0&i&0\\
0&i&0&0\\
-i&0&0&0
\end{array}
\right),
\]
where $\gamma_0\cdot\gamma_0=1$, $\gamma_i\cdot\gamma_i=-1$ for all $i=1,\ldots,4$
and $\gamma_i\cdot\gamma_j=-\gamma_i\cdot\gamma_j$ for all $i\neq j$.

The Lorentzian manifold $(\tilde{B}_a,g_a)$ with $C^1$-metric is simply connected and oriented. 
There exists a 
unique spin structure
\[
\pi: Spin(\tilde{B}_{a})\to SO(\tilde{B}_{a}),
\]
whereby $Spin(\tilde{B}_{a})$ denotes a 
$\Spin(1,4)$-principal fibre bundle over $\tilde{B}_{a}$, the spinor frame bundle, which is a $\mathbb{Z}_2$-covering 
of the orthonormal
frame bundle $SO(\tilde{B}_a)$ such that the fibre actions of $\Spin(1,4)$ and $\SO(1,4)$ are compatible  with the 
projections $\pi$ and $\lambda$. 
We denote the spinor bundle on $(\tilde{B}_a,g_a)$ by \[S:=Spin(\tilde{B}_a)\times_{\Spin(1,4)} \Delta_{1,4}\ .\] 
The spinor bundle $S$ is globally trivial on $\tilde{B}_a$.
With respect to a $C^1$-section of $\pi: Spin(\tilde{B}_a)\to \tilde{B}_a$ (i.e., a global spinor frame 
of class $C^1$) the space $C^1(\tilde{B}_a,S)$ of differentiable spinor fields is uniquely identified with the space 
$C^1(\tilde{B}_a,\Delta_{1,4})$ of 
$\Delta_{1,4}$-valued differentiable functions on $\tilde{B}_a$.
The spinor bundle $S$ admits an invariant inner product, 
which we denote by $\langle\cdot,\cdot\rangle_S$, and 
$\cdot:TM\otimes S\to S$ denotes 
the Clifford multiplication of tangent vectors with spinors.
The spinor derivative is $\nabla^S$ and
the Dirac operator $D^S$ acting on spinor fields is given with respect to some (local) orthonormal frame 
$\{t_0,\ldots,t_4\}$ by
\[\begin{array}{rcl}
D^S\ :\quad C^1(\tilde{B}_a,S)&\to&\ \ C^0(\tilde{B}_a,S)\quad .\\[2mm]
\phi\quad{ } &\mapsto& -t_0\cdot \nabla^S_{t_0}\phi\ +\ 
\sum_{i=1}^4 t_i\cdot \nabla^S_{t_i}\phi\end{array}\]

The twistor equation for spinor fields $\phi$ on $\tilde{B}_a$ is given by
\[
\nabla^S_X\phi+\frac{1}{n}X\cdot D^S\phi=0\quad\qquad\mbox{for\ all}\quad X\in T\tilde{B}_a\ .
\]
This is a first order differential equation on spinors, which is well known
to be 
conformally covariant (cf. section \ref{ab3}).
We call a spinor field $\phi\in C^1(\tilde{B}_a,S)$ a twistor spinor if it satisfies the twistor equation.
Any spinor field $\phi\neq 0$ defines in a unique way a non-trivial vector 
field $V_\phi$ (spinor square) by demanding the 
relation
\[
g_a(V_\phi,X)=\langle\phi,X\cdot\phi\rangle_S\qquad\quad \mbox{for\ all\ }\ X\in T\tilde{B}_{a}\ .
\]
If $\phi$ is a twistor spinor then $V_\phi$ is a conformal Killing vector field, i.e.,
it holds
\[L_{V_\phi}g_a=\frac{2}{n}div(V_\phi)\cdot g_a\] 
for the Lie derivative of the metric along $V_\phi$.

The $C^1$-frame $e:B^>_a\to SO(\tilde{B}_a)$ (cf. Lemma \ref{LMH}) admits exactly two lifts 
(of class $C^1$) to
the spinor frame bundle 
$Spin(\tilde{B}_a)$.
We choose one of these lifts and denote it by \[
e_s:B^>_a\to Spin(\tilde{B}_a)\ .\] 
Any spinor 
field $\phi$ on 
$B^>_a$ can then 
be uniquely represented with respect to the spinor frame $e_s$ by a $\Delta_{1,4}$-valued function $w$. 
It is
\[\phi\ =\ [\ e_s\ ,\ w\ ]\quad \in C^1(B^>_a,S) \]
for some function $w\in C^1(\tilde{B}_a,\Delta_{1,4})$.
Now let $w(b,c)$ denote the constant $\Delta_{1,4}$-valued function $(b,-c,0,0)^\top$,
where $(b,c)\in \CC^2$.
We  set
\[
\psi^>_{bc}:=(x_0e_0+re_1)\cdot[\ e_s\ ,\ w(b,c)\ ]\qquad\quad \mbox{on}\ \ B^>_a\ ,
\]
whereby the dot $\cdot$ denotes Clifford multiplication.
Obviously, the spinor field $\psi^>_{bc}$ is an element of $C^1(B^>_a,S)$ for all $(b,c)$. 
Calculating the Clifford product 
results to
\[
\psi^>_{bc}\ =\ [\ e_s\ ,\ \left(\begin{array}{c}-x_0b\ \ \  \\x_0c\\rb\\rc\end{array}\right)\ ]\ .
\]
We denote by 
\[
V\ :=\ -2x_0r\frac{\d}{\d r}\ -\ (r^2+x_0^2)\frac{\d}{\d x_0}
\]
a smooth vector field on $\RR^5$.

\begin{THEO} \label{TH1} Let $(b,c)\in\CC^2$ and $\psi^>_{bc}$ a spinor on $(B^>_a,g_a)$ with 
$a>0$.
\begin{enumerate}
\item The spinor field $\psi^>_{bc}$ on $B^>_a$ admits a unique extension $\psi_{bc}$ 
to $(\tilde{B}_a,g_a)$ of class $C^1$.
\item The unique extension $\psi_{bc}$ is a twistor spinor on $(\tilde{B}_a,g_a)$.
\item For $(b,c)\neq 0$ the twistor spinor $\psi_{bc}$ admits exactly one zero at the origin $\{0\}\in \tilde{B}_a$.
\item The zero set of the spinor length square $u_{bc}:=\langle\psi_{bc},\psi_{bc}\rangle_S$ is $L_o$. The function
$u_{bc}$ solves the equation
\[
-u_{bc}\cdot Ric^0=(n-2)\cdot Hess(u_{bc})^0
\] 
on $\tilde{B}_a\smallsetminus L_o$, 
where $Ric^0$ and $Hess(u_{bc})^0$ denote the trace-free parts of the symmetric tensors $Ric^{g_a}$ 
resp. 
$Hess^{g_a}(u_{bc})$.
In particular, the metric $\tilde{g}_a=\frac{1}{u_{bc}^2}g_a$ is Einstein for $u_{bc}\neq 0$. 
\item The spinor square $V_{\psi_{bc}}$ is a smooth conformal vector field on $(\tilde{B}_a,g_a)$.
It holds
\[
V_{\psi_{bc}}= (b^2+c^2) \cdot V\ .
\]
\item
The vector $V_{\psi_{bc}}$ is timelike on $\tilde{B}_a\smallsetminus L_o$, 
lightlike on $L_o\smallsetminus\{0\}$ and zero only in the origin 
$\{0\}\in \tilde{B}_a$.
\end{enumerate}
\end{THEO} 
Here a vector
$X\neq 0$ on $(\tilde{B}_a,g_a)$ is
called timelike if $g_a(X,X)<0$ and lightlike if $g_a(X,X)=0$.
In short, Theorem \ref{TH1} says that there exists a $2$-dimensional set of twistor spinors on $(\tilde{B}_a,g_a)$
for all $a>0$, which admit an isolated zero at the origin. There exist no further twistor spinors on
$(\tilde{B}_a,g_a)$, since the Eguchi-Hanson metric admits exactly two linearly independent (parallel) twistor spinors for $a>0$. 

For $a=0$ we set $\tilde{B}_a=\RR^5$ and $g_a=g_0$. All twistors 
with zero at the origin on $(\RR^5,g_0)$ are given by
\[ 
\psi_{w_0}\ =\ \left(\sum_{i=0}^4 x_i\frac{\d}{\d x_i}\right)\cdot [\ u_s\ , \ w_0\ ]\ ,
\] 
where $u_s$ is a lift of the standard frame $\{\frac{\d}{\d x_0},\ldots,\frac{\d}{\d x_4}\}$ 
to $Spin(\tilde{B}_a)$
and $w_{0}\in\CC^4$ a constant (cf. \cite{BFGK91}). The metric $g_0$ is smooth and flat on $\tilde{B}_0$.
For $a > 0$ the situation changes. 
In this case the metric $g_a$ is only of class $C^1$ and the space of twistors shrinks to dimension $2$. 
From Proposition \ref{PR2} we know that $g_a$ with twistors 
$\psi_{bc}$ gives rise to a curved (conformal) geometry on $\tilde{B}_a$.
In fact, $g_a$ is not conformally flat in any neighbourhood
of the origin $\{0\}\in \tilde{B}_{a}$. Nevertheless, the twistors $\psi_{bc}$ have a zero at the origin.  
This is an important observation for our construction.

\begin{CO} \label{CO1} There exists a family of Lorentzian $C^1$-metrics $g_a$, $a>0$, in dimension $5$, which admit 
twistor
spinors and a smooth causal conformal Killing vector field, all with isolated zero at some 
point $\{p\}$ such that 
$g_a$ is 
non-conformally flat around the zero at $\{p\}$.
\end{CO}

We remark that
the vector field $V$ is complete on $\tilde{B}_a$, i.e., the flow of $V$
to the time $t$ generates a $1$-parameter group of conformal transformations on $\tilde{B}_a$.
These conformal transformations are not isometries with respect to any metric in the conformal
class $c_a:=[g_a]$ (cf. section \ref{ab3}). Conformal transformations with the latter property are called essential. 
In particular, the conformal Killing vector field $V$ is called essential.
The statement of Corollary \ref{CO1} implies the existence of essential conformal Killing 
fields and transformations on non-compact Lorentzian spaces, which are not
conformally flat.
A conjecture by A. Lichnerowicz states that essential conformal transformation groups do not
exist on any compact Lorentzian manifold unless it is conformally flat (cf. \cite{DAG91}). 
In fact, we do not expect that our construction works on compact spaces.

We want to add some further comments concerning our construction. 
The metric $g_a$ can be considered as a completion of the metric $-ds^2+g_{EH}$, 
which is Ricci-flat and 'asymptotically Minkowskian', to the set $L$ with infinity $L_o$.
The twistors extend to $L$ as well with a zero at some point of infinity.
In general, it is known that a Lorentzian metric with differentiable Weyl tensor 
has to be conformally flat in the causal 
past and future 
of a zero of a twistor spinor (cf. section \ref{ab3}).
Therefore, it is also reasonable in our construction to do the conformal
completion to $L$ by using the flat Minkowski metric $g_0$ on the 'other side' of the infinity set $L_o$. 
There exists no 
extension (conformal completion) with differentiable Weyl tensor 
of $g_a$ on $B_a$ to a neighbourhood of the origin, which 
is not conformally flat
on $L$, but preserves the existence
of a twistor spinor.
This fact implies that our completion of $g_a$ can not be analytic.
We want to point out again that our construction is even not of class $C^\infty$.
However, it remains the question whether there is a conformally equivalent metric to $g_a$ on $\tilde{B}_a$,
whose regularity is better then of class $C^1$.
The existence of a $C^1$-extension of the Weyl tensor of $g_a$ to the infinity set $L_o$
certainly does not pose an obstruction to this question.

\section{Proof of statements}
\label{ab3}

We prove here the statements which we made in the previous section.
We start with a discussion of the differentiability of certain functions on $\tilde{B}_a\subset \RR^5$.
For some arbitrary $p$-tuple $I_p=(i_1,\ldots,i_p)\in \{0,\ldots,4\}^p$ let us denote by
\[
\partial_{I_p}:=\frac{\d}{\d x_{i_1}}\cdots\frac{\d}{\d x_{i_p}}
\]
a partial derivative of order $p$. Moreover, for any $5$-tuple $l=(l_r,l_0,\ldots,l_4)$ with 
$l_r,l_0,\ldots,l_4\in \mathbb{N}\cup\{0\}$
we set $s_l:=-l_r+\sum_{i=0}^4l_i$ and define the smooth function
\[
f_l=f(l_r,l_0,\cdots, l_4):= r^{-l_r}\cdot x_0^{l_0}\cdot\ldots \cdot x_4^{l_4}\qquad\qquad \mbox{on}\ \ B_a\ .
\]
We say that the rational function $f_l$ is of order $s_l$.
Remember that we defined the radial function $r_o$ to be $(r^2-x_0^2)/r$ on $B_a$ and identically 
zero on $L$ (cf. section \ref{ab2}). 
For any function $f$ on $B_a$ we understand the product 
$r_o\cdot f$ in a unique way as a function on $\tilde{B}_a=B_a\cup L$, 
which is identically zero on $L$. For $t>0$ a real number we denote \[B_a^t := B_a\cap \{x\in\RR^5:\ r\leq t\}\ .\]
Notice that if a function $f$ is continuous on $B_a$ and its absolute value $|f|$ is bounded on $B_a^t$ for all $t>0$ then 
$r_o\cdot f$
is continuous on $\tilde{B}_a$. In fact, for this conclusion it is sufficient for 
$|f|$ to be bounded on $B_{\tilde{a}}^t$ for all $t>0$ with some $\tilde{a}>a$.

\begin{LM} \label{LM1} Any function on $\tilde{B}_a$ of the form $r_o^m\cdot f_l$ 
with $m>0$ is of class $C^{k-1}$ but not of class $C^{k}$, where $k:=\mathrm{min}\{m,m+s_l\}$. 
\end{LM}

{\bf Proof.} First, we note that $|x_i|<r$ on $B_a$ for all $i=0,\ldots,4$, and we see that
the absolute value $|f_l|$ of any function of the form $f_l$ with $s_l\geq 0$ is bounded on $B_a^t$ by $t^{s_l}$ 
for all $t>0$.
More generally, the absolute value of the partial derivative $\d_{I_p} f_l$ is bounded on $B_a^t$ for all $t>0$ 
if $s_l-p\geq 0$. 
In particular, the absolute value of $x_0/r$ is bounded on $B_a$. Moreover, $x_0/r$ is continuous on $(B_a\cup L_o)\smallsetminus\{0\}$.
This shows that the extension of the function $r-x_0\cdot x_0/r$\ by zero to the origin in $\RR^5$  
is a continuous function on $B_a\cup L_o$ and is identically zero on $L_o$. And this implies  
that the coordinate $r_o$ is continuously defined on $\tilde{B}_a$. 
The function $r_o$ is not
continuously differentiable. However, it is
\[
dr_o=\frac{-2x_0}{r}dx_0\ +\ \sum_{i=1}^4 \left(\frac{x_ix_0^2}{r^3}+\frac{x_i}{r}\right)dx_i\ ,
\]
and we see that the coefficients of $dr_o$ are bounded on $B_a$.
(This implies that $d(r_o^2)=2r_odr_o$ is continuous on $\tilde{B}_a$, i.e., $r_o^2$ is of class $C^1$.) 
  
Eventually, since a function of the form $r_o^m\cdot f_l$ admits terms of lowest
order $m+s_l$, such a function is at most of class $C^{m+s_l-1}$. However, a $p$th order derivative of 
$r_o^m\cdot f_l$ on $B_a$ can be extended continuously
to the zero function on $L$ only if $p<m$. In fact, any application of a derivative $\d_{I_m}$ of order $m$ admits 
a non-trivial
term of the form $m!\cdot f_l\cdot \Pi_{k=1}^m dr_o(\frac{\d}{\d x_{i_k}})$, which can not be extended continuously
by zero to $L$. On the other side, it is $\d_{I_{m-1}}(r_o^m\cdot f_l)=r_o\cdot h$, whereby $|h|$ is bounded on $B_a^t$ 
for all $t>0$ if $s_l\geq 0$. 
This shows that $r_o^m\cdot f_l$ is of class $C^{k-1}$ with $k:=\mathrm{min}\{m,m+s_l\}$, but it is not $k$-times
continuously differentiable.	\hfill$\Box$\\

Now we set \[
\omega_a:= (ar_o)^4\cdot (r\sigma_3)^2\qquad\quad\mbox{and}\qquad\quad \rho_a=\frac{a^4r_o^2}{1-(ar_o)^4}\cdot 
\left(\frac{\alpha}{r}\right)^2\ .
\]
With these notations it is
$g_a=g_0-\omega_a+\rho_a$ on $\tilde{B}_a$, where $g_0$ is the flat Minkowski metric on $\tilde{B}_a$.\\

{\bf Proof of Proposition \ref{PR1}.} The metric $g_0$ is smooth on $\tilde{B}_a$. We have to discuss
the differentiability of $\omega_a$ and $\rho_a$. The coefficients of the $1$-forms
$r\cdot \sigma_3$ and $\alpha/r$ are of order $s_l=0$ resp. $s_l=1$. 
With 
application
of Lemma \ref{LM1} we conclude that  $\omega_a$ is of class $C^3$ and $\rho_a$ is of class $C^1$.
The symmetric $2$-form $\rho_a$ is not of class $C^2$. This implies that the symmetric bilinear form
$g_a$ on $\tilde{B}_a$ is of class $C^1$ for all $a>0$, but it is not of class $C^2$.

We postpone the proof that $g_a$ is a metric of Lorentzian signature until the proof of Proposition \ref{PR2}.
The proof of Lemma \ref{LMH} about the existence of the orthonormal frame $e$ will show the Lorentzian signature of $g_a$ as well.
\hfill $\Box$\\

For the proof of Proposition \ref{PR2} we use the coordinate change
\[ 
\begin{array}{lccl}
\Psi:&\RR^5\smallsetminus L_o &\quad\to\quad&\quad \RR^5\smallsetminus L_o\ ,\\
&(x_0,r,\varphi_i)&\quad\mapsto\quad&\quad (s,R,\varphi_i)=\left(\ \frac{-x_0}{r^2-x_0^2}\ ,\ \frac{r}{r^2-x_0^2}\ ,\
\varphi_i\ \right)\end{array}
\]
whereby the $\varphi_i$'s are some (local) coordinates on $S^3$ which remain unchanged. The coordinate 
transformation  $\Psi$ 
is smooth on $\RR^5\smallsetminus L_o$ and it holds
\[
\begin{array}{ccrcr}
dx_0&=& -\frac{s^2+R^2}{\ (R^2-s^2)^2}ds&+& \frac{2sR}{\ (R^2-s^2)^2} dR,\\[4mm]
dr&=&\ \  \frac{2sR}{\ (R^2-s^2)^2}ds&-&  \frac{s^2+R^2}{\ (R^2-s^2)^2} dR,\\[4mm]
\frac{\d}{\d r}&=& -2sR\frac{\d}{\d s}&-& (s^2+R^2)\frac{\d}{\d R},\\[4mm]
\frac{\d}{\d x_0}&=& -(s^2+R^2)\frac{\d}{\d s}&-& 2sR\frac{\d}{\d R}\ .
\end{array}
\]
This shows also $T=\frac{\d}{\d R}$ and $V=\frac{\d}{\d s}$.\\

{\bf Proof of Proposition \ref{PR2}:} First, we calculate the symmetric bilinear form $g_a$ on $B_a$ with respect to the coordinate 
transformation $\Psi$. 
Remember that $\alpha=(x_0^2+r^2)dr-2x_0rdx_0$. It holds
\[\begin{array}{ccl}
\alpha &=& \frac{-dR}{\ (R^2-s^2)^2}\ ,\\[2mm]
-dx_0^2+dr^2 &=& \frac{-ds^2+dR^2}{\ (R^2-s^2)^2}\ . 
\end{array}
\]
With $R^2=r_o^{-2}$ on $B_a$ and $r^2-x_0^2=(R^2-s^2)^{-1}$ we obtain
\begin{eqnarray*}
g_a&=&\ \ \frac{1}{\ (R^2-s^2)^2}(\ -ds^2+dR^2+R^2(\sigma_1^2+\sigma_2^2+(1-(a/R)^4)\sigma_3^2)\ )\\
&&- \frac{R^2\cdot dR^2}{r^2(1-(R/a)^4)\cdot(R^2-s^2)^4}\\[4mm]
&=&\frac{1}{\ (R^2-s^2)^2}\left(\ -ds^2+dR^2+R^2(\sigma_1^2+\sigma_2^2+(1-(a/R)^4)\sigma_3^2)\right.\\
&&\qquad\qquad\qquad- \frac{1}{1-(R/a)^4}dR^2\ \big)\\[4mm]
&=& \frac{1}{\ (R^2-s^2)^2}\big(\ -ds^2+\frac{dR^2}{1-(a/R)^4}\ +\ R^2(\sigma_1^2+\sigma_2^2+(1-(a/R)^4)\sigma_3^2\ \big)
\end{eqnarray*}
and we can conclude that
\[
\tilde{g}_a\ =\ \frac{1}{(r^2-x_0^2)^2}\cdot g_a\ =\ -ds^2\ +\ g_{EH}
\]
on $B_a$. The corresponding (even simpler)
calculation on $L\smallsetminus L_o$, where $r_o\equiv 0$, 
shows that
\[  
\tilde{g}_a\ =\ \frac{1}{(r^2-x_0^2)^2}g_a\ =\ -ds^2\ +\ 
dR^2\ +\ R^2(\sigma_1^2+\sigma_2^2+\sigma_3^2)\ ,
\]
i.e., $\tilde{g}_a$ is the flat metric for $r<|x_0|$.
In particular, since $\tilde{g}_a$ on $\tilde{B}_a\smallsetminus L_o$ is a metric of Lorentzian signature, we have shown that the 
conformally equivalent symmetric bilinear form $g_a$ of class $C^1$ on 
$\tilde{B}_a$ is a metric and admits
Lorentzian signature as well, which completes the proof of Proposition \ref{PR1}.

As next we review curvature properties of the Eguchi-Hanson metric $g_{EH}$. This discussion 
will provide us with all the information that we need to prove our claims about the curvature properties of 
the Lorentzian metrics $\tilde{g}_a$ and $g_a$. 
Let us fix the orthonormal frame
\[
\{f_1,f_2,f_3,f_4\}:=\left\{\ -\beta\frac{\d}{\d R}\ ,\ R^{-1}\frac{\d}{\sigma_1}\ ,\ R^{-1}\frac{\d}{\sigma_2}\ ,\
(R\beta)^{-1}\frac{\d}{\sigma_3}\ \right\}\ ,
\]
where $\beta:=\sqrt{1-(a/R)^4}$. We denote by $\{f^i:i=1,\ldots,4\}$ the dual frame. The
connection $1$-form $\omega$ and the curvature $2$-form $\Omega$ of the Levi-Civita connection
$\nabla^{g_{EH}}$ are determined by the structure equations
\[
df^i=\sum_{k=1}^4\omega^i_k\wedge f^k\qquad \mbox{and}\qquad
\Omega_j^i=d\omega^i_j-\sum_{k=1}^4 \omega^i_k\wedge \omega^k_j\ .
\]
It holds 
\[
\omega^i_j=g_{EH}(\nabla^{EH}f_i,f_j)\qquad \mbox{and}\qquad
\Omega_j^i=g_{EH}(R(e_i,e_j)\cdot,\cdot)\ ,
\]
whereby 
\[R(e_i,e_j)=\nabla^{EH}_{e_i}\nabla^{EH}_{e_j}-\nabla^{EH}_{e_j}\nabla^{EH}_{e_i}-\nabla^{EH}_{[e_i,e_j]}\ .\]
The components are explicitly calculated as
\[
\begin{array}{rcrcccc}
\omega^1_2&=&\omega^3_4&=&-\beta R^{-1}\cdot f^2&=&-\beta\cdot \sigma^1,\\[3mm]
\omega^1_3&=&-\omega^2_4&=&-\beta R^{-1}\cdot f^3&=&-\beta\cdot \sigma^2,\\[3mm]
\omega^1_4&=&\omega^2_3&=&-\gamma\cdot f^4&=&-\gamma R\beta\cdot \sigma^3,
\end{array}
\]
whereby $\gamma=\beta R^{-1}+\beta'$ and $\beta'=\frac{\d \beta}{\d R}=2(a/R)^4(R\beta)^{-1}$, and 
\[
\begin{array}{l}
\Omega^1_2=\ \ \Omega^3_4=-\frac{2a^4}{R^6}\lambda^1_-\ ,\\[3mm]
\Omega^1_3=-\Omega^2_4=-\frac{2a^4}{R^6}\lambda^2_-\ ,\\[3mm]
\Omega^1_4=\ \ \Omega^2_3=\quad \frac{4a^4}{R^6}\lambda^3_-\ ,
\end{array}
\]
where the $\lambda^i_-$'s build a basis of the anti-selfdual $2$-forms for $g_{EH}$ and are defined as
\[
\begin{array}{l}
\lambda^1_-=f^1\wedge f^2-f^3\wedge f^4\ ,\\[2mm]
\lambda^2_-=f^1\wedge f^3-f^4\wedge f^2\ ,\\[2mm]
\lambda^3_-=f^1\wedge f^4-f^2\wedge f^3\ .
\end{array}
\]
It follows that the Riemannian curvature tensor $R^{EH}$ of $g_{EH}$ is anti-selfdual.
This implies that $g_{EH}$ is Ricci-flat and $R^{EH}$
equals the Weyl tensor $W^{EH}$, i.e., we have
\[ R^{EH}=W^{EH}=W^{-}\neq 0\ .
\] 
In particular, since the Weyl tensor is a complete obstruction
to conformal flatness in dimension $4$, we can see that $g_{EH}$ is nowhere conformally flat on its domain of
definition (which is $B_a\cap E$ resp. $\Psi(B_a\cap E)$\ ).   

Now the metric $\tilde{g}_a=-ds^2+g_{EH}$ is an ordinary semi-Riemannian product. Hence the curvature components
of $\tilde{g}_a$ in direction  of the coordinate $\frac{\d}{\d s}$ vanish, i.e.,
the curvature tensor of $-ds^2+g_{EH}$ is entirely determined by
the components of the Riemannian curvature tensor $R^{EH}$.
In particular, we see that the metric  $-ds^2+g_{EH}$ is Ricci-flat and the components of the Weyl tensor
$W^{\tilde{g}_a}$ of $\tilde{g}_a$ in direction of the coordinate $\frac{\d}{\d s}$ do vanish as well. 
Since, by construction, the metric $\tilde{g}_a=-ds^2+g_{EH}$ is conformally equivalent to $g_a$ on $B_a$,
we know yet the Weyl tensor $W^{g_a}$ of $g_a$ on $B_a$ as well. It is simply a rescaling of $W^{\tilde{g}_a}$. 
Obviously, the metric $g_a$ 
is not conformally flat on $B_a$.
On $L\smallsetminus L_o$ the metric $g_a$ is flat and therefore conformally flat, i.e., $W^{g_a}\neq 0$ on $B_a$ and
$W^{g_a}\equiv 0$ on $L\smallsetminus L_o$.

Finally, on the lightcone $L_o$ the Weyl tensor of $g_a$ is not defined in the usual way, because $g_a$ is
only of class $C^1$ at $L_o$. 
We aim to show that the Weyl tensor of $g_a$ on $\tilde{B}_a\smallsetminus L_o$
admits a continuous extension to $L_o$. For this we note that the Weyl tensor rescales explicitly by 
$W^{g_a}=r_o^4r^4\cdot W^{\tilde{g}_a}$. Then calculating the components of
$W^{g_a}$ with respect to the coordinate system $u:=\{\frac{\d}{\d x_0},\ldots,\frac{\d}{\d x_4}\}$
using our formulae for $W^{EH}$ from above results to expressions of the form
\[
W^{g_a}\Big(\frac{\d}{\d x_i},\frac{\d}{\d x_j},\frac{\d}{\d x_k},\frac{\d}{\d x_l}\Big)
= \left\{\begin{array}{cl}A\cdot r_o^2\ +\ B\cdot r_o^6/r^4&\quad \mbox{on}\ \ B_a\\[3mm]
0&\quad \mbox{on}\ \ L\smallsetminus L_o\end{array}\right. 
\]
for all $i,j,k,l\in \{0,\ldots,4\}$, where $A,B$ are sums of functions of the form $f_l$, $\beta\cdot f_l$
and $\beta^{-1}\cdot f_l$ with order $s_l=4$,
i.e., the extensions of all components to $L_o$ by zero are $C^1$-functions on $\tilde{B}_a$.
We conclude that the Weyl tensor $W^{g_a}$ has a continuous extension of class $C^1$ on $\tilde{B}_a$. 
\hfill $\Box$\\

Now we consider the frame $e=\{e_0,\ldots,e_4\}$, which we have defined in section \ref{ab2} and which was claimed 
there to 
be orthonormal for $g_a$ on $B^>_a$ and of class $C^1$.\\

{\bf Proof of Lemma  \ref{LMH}:} 
First, we show that the frame $e$ is orthonormal in every point of $(B^>_a,g_a)$.
Obviously, this is true on $L\smallsetminus \{r=0\}$, since $g_a$ is the flat Minkowski metric
thereon. It is also obvious that the vectors $e_2,e_3$ and $e_4$ are orthonormal for $g_a$ on $B_a$
and that they are orthogonal to the remaining basis vectors $e_0$ and $e_1$.
For the latter we find with $\frac{a^4r_o^2}{1+\beta}=R^2(1-\beta)$ and $T=\frac{\d}{\d R}$ the expressions
\[
\begin{array}{lrlcl}
e_0&=&-(s^2+R^2)\frac{\d}{\d s}&-&2sR\beta\frac{\d}{\d R}\quad\qquad\mbox{and}\\[3mm]
e_1&=&-2sR\frac{\d}{\d s}&-&(s^2+R^2)\beta\frac{\d}{\d R}\ ,
\end{array}
\]
from which we see that $e_0$ and $e_1$ are orthonormal 
with respect to $g_a = (R^2-s^2)^{-2}(-ds^2+g_{EH})$ on $B_a$ as well.
We conclude that the frame $e$ is a pointwise orthonormal basis on $B_a$.

It remains to discuss the differentiability of the coefficients of the vectors $\{e_0,\ldots, e_4\}$.
For this we notice that the function
$\frac{a^4r_o^2}{1+\beta}$
is only of class $C^1$ on $B^>_a$. The function $\beta^{-1}$ is of class $C^3$ and all other functions,
which are involved in the coefficients are smooth on $B^>_a$.\hfill$\Box$\\

Let us introduce the vectors 
\[\begin{array}{l}
\tilde{e}_0:=\frac{-1}{R^2-s^2}(\ (S^2+R^2)\frac{\d}{\d s}\ +\ 2sR\beta\frac{\d}{\d R}\ )\qquad
\mbox{and}\\[3mm]
\tilde{e}_1:=\frac{-1}{R^2-s^2}(\ 2sR\frac{\d}{\d s}\ +\ (s^2+R^2)\beta\frac{\d}{\d R}\ )
\end{array}\]
with respect to the $\Psi$-transformed coordinates, and let us denote
\[
\tilde{e}:=\left\{\ \tilde{e}_0\ ,\ \tilde{e}_1\ ,\ R^{-1}\cdot\frac{\d}{\d \sigma_1}\ ,\ R^{-1}\cdot
\frac{\d}{\d \sigma_2}\ 
, \ 
(R\beta)^{-1}\cdot\frac{\d}{\d \sigma_3}\ \right\}\ ,
\]
which is an orthonormal frame with respect to $\tilde{g}_a=\frac{1}{(r^2-x_0^2)^2}\cdot 
g_a$ on $B^>_a\smallsetminus L_o$. 
As we know from the proof of Lemma \ref{LMH}, it holds $\Psi_*(e_i)=(R^2-s^2)\cdot\tilde{e}_i$, $i=0,\ldots,4$. 
Moreover, we set
\[
f:=\left\{\ -\frac{\d}{\d s}\ ,\ -\beta\frac{\d}{\d R}\ ,\ R^{-1}\cdot\frac{\d}{\d \sigma_1}\ ,\ 
R^{-1}\cdot\frac{\d}{\d \sigma_2}\ , \
(R\beta)^{-1}\cdot\frac{\d}{\d \sigma_3}\ \right\}
\] 
on $B^>_a\smallsetminus L_o$. On $B_a$ this is just the extension by $f_0$ of the frame $\{f_1,f_2,f_3,f_4\}$ that we introduced already
for the Eguchi-Hanson metric $g_{EH}$. 
The frames $\tilde{e},f$ are transformed on $B^>_a\smallsetminus L_o$ by the matrix 
\[
\kappa=\frac{1}{R^2-s^2}\left(\begin{array}{ccccc}
s^2+R^2 & 2sR & 0 & 0 & 0 \\
2sR & s^2+R^2 & 0 & 0 & 0 \\
0 & 0 &  R^2-s^2 & 0 & 0 \\
0 & 0 & 0 & R^2-s^2 & 0 \\
0 & 0 & 0 & 0 & R^2-s^2
\end{array}\right),
\]
namely it holds $\tilde{e}=f\cdot \kappa$.
With $t:=\ln \frac{R-s}{R+s}$ and \[E_{01}:=\left(\begin{array}{ccccc}
0& -1 & 0 & 0 & 0 \\
-1 & 0 & 0 & 0 & 0 \\
0 & 0 & 0 & 0 & 0 \\
0 & 0 & 0 & 0 & 0 \\
0 & 0 & 0 & 0 & 0
\end{array}\right)\] we have $\kappa=\exp(tE_{01})$. The elements in the preimage of $\kappa$ by the group 
covering
$\lambda:\Spin(1,4)\to \SO(1,4)$  
are given by  $\pm \exp(\frac{t}{2}\gamma_0\gamma_1)$, whereby we use the $\gamma$-matrices introduced in section 
\ref{ab2}. We choose in the following $\tilde{\kappa}:=\exp(\frac{t}{2}\gamma_0\gamma_1)$, which is given by 
\[
\tilde{\kappa}=\frac{1}{\sqrt{R^2-s^2}}\left(\begin{array}{cccc}
R& 0 & -s & 0 \\
0 & R & 0 & s\\
-s & 0 & R & 0 \\
0 & s & 0 & R 
\end{array}\right).
\] 

Before we start with the proof of Theorem \ref{TH1}, let us recall the 
conformal covariance of the twistor equation in explicit terms. In general, let $\tilde{g}=e^{2\sigma}g$ be a rescaled
metric in the conformal class of a given metric $g$ and let $\varphi=[v_s,w]$ be a twistor with respect 
to $g$, whereby  
$v_s$ denotes the lift of some orthonormal frame $v$. We set $\tilde{\varphi}:=[\tilde{v}_s,w]$, 
whereby $\tilde{v}_s$ denotes the lift of the
rescaled frame $\tilde{v}:=v\cdot (e^{-\sigma}id)$, which naturally corresponds to the lift $v_s$.
Then the spinor field $e^{\sigma/2}\cdot \tilde{\varphi}$ is a twistor spinor with respect to the 
rescaled metric $\tilde{g}$
(cf. \cite{BFGK91}).\\

{\bf Proof of Theorem \ref{TH1}:} The verification of the first two statements of Theorem \ref{TH1}
is the main part of the proof. We  
will show this in some few steps. 
First, we prove that
$\psi^>_{bc}$ is a twistor on $B_a$ and also on $L\smallsetminus (L_o\cup \{r=0\})$, which already implies that $\psi^>_{bc}$
is a twistor on $B^>_a$. Thereby, we will not directly check the twistor equation
for $\psi^>_{bc}$, but first
use the conformal transformation from $g_a$ to the Ricci-flat metric $\tilde{g}_a$. 
In the next step we
show that 
$\psi^>_{bc}$ extends to a $C^1$-spinor on $\tilde{B}_a\smallsetminus\{0\}$. 
This spinor will still be 
a twistor. Finally, we show that the latter spinor can be extended to the origin by a zero. The resulting 
spinor $\psi_{bc}$ is a unique continuous extension of $\psi^>_{bc}$, which is of class $C^1$ and solves the twistor 
equation 
everywhere on $\tilde{B}_a$.

To start with, let us consider $\psi^>_{bc}$ on $B^>_a\smallsetminus L_o$. The spinor $\psi^>_{bc}$ is given 
with respect to the spinor frame $e_s$ by $[\ e_s\ ,\ (-x_0b,x_0c,rb,rc)^\top\ ]$.
It holds $e=\tilde{e}\cdot((R^2-s^2)id)$, where $\tilde{e}$ is orthonormal with respect to
$\tilde{g}_a=\frac{1}{(r^2-x_0^2)^2}g_a$. Let $\tilde{e}_s$ be the corresponding lift of the rescaled frame 
$\tilde{e}$.
Then the spinor 
\[
\nu_{bc}:=\sqrt{R^2-s^2}\cdot\tilde{\psi}^>_{bc}=[\ \tilde{e}_s\ ,\ \sqrt{R^2-s^2}\cdot (-x_0b,x_0c,rb,rc)^\top\ ]
\]	
is a twistor with respect to $\tilde{g}_a$ (by conformal covariance). Further, it holds
\[
\nu_{bc}=\sqrt{R^2-s^2}\cdot [\ f_s\ ,\ \tilde{\kappa}\cdot (-x_0b,x_0c,rb,rc)^\top\ ]\ , 
\]
where $f_s$ denotes the lift of the frame $f$, which corresponds to the lift $\tilde{e}_s$.
Eventually, with \[
\sqrt{R^2-s^2}\cdot\tilde{\kappa}\left(\begin{array}{c}-x_0b\ \ \ \\x_0c\\rb\\rc\end{array}\right)=
\left(\begin{array}{cccc}
R& 0 & -s & 0 \\
0 & R & 0 & s\\
-s & 0 & R & 0 \\
0 & s & 0 & R
\end{array}\right)\left(\begin{array}{c}-x_0b\ \ \ \\x_0c\\rb\\rc\end{array}\right)
=\left(\begin{array}{c}0\\0\\b\\c\end{array}\right)
\]
we find that
\[
\nu_{bc}=[\ f_s\ ,\ (0,0,b,c)^\top\ ]\ .
\]

The spinor derivative of $\nu_{bc}$ with respect to $\tilde{g}_a$ is given by
\[
\tilde{\nabla}^S\nu_{bc}=[\ f_s\ ,\ \frac{1}{2}\cdot \sum_{0\leq i<j\leq 4}\omega^i_j\otimes\gamma_i\gamma_j\cdot 
(0,0,b,c)^\top\ ]\ ,
\]
where the $\omega^i_j$'s are the components of the Levi-Civita connection of $\tilde{g}_a$.
On $B_a$ it holds $\omega^0_j=0$ and the other $\omega^i_j$'s are just the components
that we calculated in the proof of Proposition \ref{PR2} for the Eguchi-Hanson metric $g_{EH}$.
Notice also that on $L\smallsetminus (L_o\cup \{r=0\})$ the components $\omega^i_j$ admit the same expressions (with $\beta\equiv 1$) 
as on $B_a$
with respect to the frame $f$.
The relations for the $\omega^i_j$'s immediately prove that
$\nu_{bc}$ is a parallel spinor with respect to $\tilde{g}_a$ on $B_a^>\smallsetminus L_o$ 
for any $(b,c)\in \CC^2\smallsetminus {0}$. (In fact, the 
spinors of the 
form $\nu_{bc}$ restricted to the Eguchi-Hanson metric $g_{EH}$, which is a hyperk{\"a}hler metric for any $a>0$, 
form the space of all parallel spinors thereon.) Any parallel spinor is a twistor spinor. In particular,
$\nu_{bc}$ is a twistor spinor for $\tilde{g}_a$. Hence, by conformal covariance and the fact that $\psi^>_{bc}$ is of class 
$C^1$ on $B_a^>$, it follows that $\psi^>_{bc}$ is a twistor on $(B_a^>,g_a)$.

Now let 
\[
G=r^{-1}\cdot \left(\begin{array}{rrrrr}
r&0&0&0&0\\
0&x_1&x_2&x_3&x_4\\
0&-x_2&x_1&-x_4&x_3\\
0&-x_3&x_4&x_1&-x_2\\
0&-x_4&-x_3&x_2&x_1\end{array}\right)
\]
be a matrix valued function on $B_a^>$.
It holds $e\cdot G=\{\frac{\d}{\d x_0},\ldots,\frac{\d}{\d x_4}\}$ on $L\smallsetminus(L_o\cup \{r=0\})$.
The standard frame $u$ is orthonormal on $L\smallsetminus(L_o\cup \{r=0\})$ and admits a smooth 
extension to $L\smallsetminus L_o$. Of course, the matrix $G$ is singular for $r=0$. A transformation matrix
for corresponding spinor frames is given by
\[
\tilde{G}=r^{-1}\cdot \left(\begin{array}{cccc}
r&0&0&0\\
0&r&0&0\\
0&0&\ \ x_1+ix_2&x_3+ix_4\\
0&0&-x_3+ix_4&x_1-ix_2
\end{array}\right)\ .
\]
This form of the matrix is due to the fact that $\Spin(4)$ is isomorphic to $\SU(2)\times \SU(2)$.
The spinor $\psi_{bc}^>=[e_s,(-x_0b,x_0c,rb,rc)^\top)]$ is presented with respect to
the spinor frame $u_s$ on $L\smallsetminus(L_o\cup \{r=0\})$ by
\[
\psi_{bc}^>=[\ u_s\ ,\ \tilde{G}^{-1}(-x_0b,x_0c,rb,rc)^\top)\ ]\ .
\]
Obviously, the vector valued function 
\[
\tilde{G}^{-1}\left(\begin{array}{c}\!\!\!\!-x_0b\\x_0c\\rb\\rc\end{array}\right)=\left(\begin{array}{c}
\!\!\!\!-x_0b\\x_0c\\ (x_1-ix_2)b-(x_3+ix_4)c\\
(x_3-ix_4)b+(x_1+ix_2)c
\end{array}\right)
\] 
is non-singular and smooth
on $L\smallsetminus L_o$. Hence the spinor $\psi_{bc}^>$ on $B_a^>$ admits a $C^1$-extension to $\tilde{B}_a\smallsetminus \{0\}$.
We denote this extension by $\psi_{bc}^o$, which is by continuity reasons a twistor on $\tilde{B}_a\smallsetminus \{0\}$.

We still have to show that $\psi_{bc}^o$ extends further to a $C^1$-spinor $\psi_{bc}$ on $\tilde{B}_a$.
For this purpose,
we improve our change of frame from above and introduce a non-singular $C^1$-frame around the origin.
Then we show that the components of $\psi^o_{bc}$ 
with respect to a corresponding non-singular spinor frame are of class $C^1$. So let 
\[
Q=\left(\begin{array}{ccccc}
k&q&0&0&0\\
q&k&0&0&0\\
0&0&1&0&0\\
0&0&0&1&0\\
0&0&0&0&1
\end{array}\right)
\]
with 
\[\begin{array}{l}
k\ :=\ \sqrt{\frac{1+r_o^2\rho}{1-4x_0^2\rho}}\ \cdot\ \left(\ 1\ -\ \frac{(r^2+x_0^2)\beta^2}{r^2(1+\beta)}\cdot\rho\ \right)     \qquad 
\mbox{and}\\[4mm]
q\ :=\ \sqrt{\frac{1+r_o^2\rho}{1-4x_0^2\rho}}\ \cdot\ \frac{2x_0(r^2+x_0^2)\beta^2}{r(1+\beta)}\ \cdot\ \rho,
\end{array}
\]
whereby $\rho=a^4\beta^{-2}r_o^2$.
It is $Q\equiv \mathbb{1}$ on $L$ and  
$4x_0^2\rho<1$ on an open neighbourhood of $L$. For $4x_0^2\rho<1$ the function $k$ is well defined and of class $C^1$ (cf. Lemma \ref{LM1}).
It follows that on a certain open neighbourhood of $L$ in $\tilde{B}_a$ the function $k$ is positive.
We denote this set by $C_a$. In fact, it holds $k^2-q^2\equiv 1$ on $C_a$, i.e.,   
the transformation matrix $Q$ takes values in $\SO_o(1,4)$ and  
is of class $C^1$ on $C_a$. The matrix $Q$ is useful, because  
the transformed frame
$\{h_0,\ldots,h_4\}:=e\cdot Q$ is given by
\[\begin{array}{l}h_0=(1-4x_0^2\rho)^{-1/2}\cdot\frac{\d}{\d x_0}\ ,\\[2mm]
h_1=(1+\rho_1)^{-1/2}\cdot\frac{\d}{\d r}\ +\ 
\rho_2(1+\rho_1)^{-1/2}\cdot(1-4x_0^2\rho)^{-1}\cdot\frac{\d}{\d x_0}\ ,\\[2mm]
h_i=e_i\qquad \mbox{for}\ \ i=2,3,4\ \end{array}\] 
on $B_a^>\cap C_a$, i.e., the first basis vector $h_0$ admits now a continuous extension to $\{r=0\}$. 
The remaining basis vectors are still singular at $\{r=0\}$. However,  
a straightforward calculation shows that the frame 
$\tilde{h}=\{\tilde{h}_0,\ldots,\tilde{h}_4\}:=e\cdot (QG)$
admits a $C^1$-extension to $\{r=0\}$, i.e., $\tilde{h}$ is a non-singular $C^1$-frame on $C_a$, which is an open 
neighbourhood of the origin. 

A corresponding transformation matrix to $Q$ for spinor frames is given by
\[
\tilde{Q}=\left(\begin{array}{cccc}
\sqrt{\frac{k+1}{2}}&0&\frac{-q}{\sqrt{2(k+1)}}&0\\
0&\sqrt{\frac{k+1}{2}}&0&\frac{q}{\sqrt{2(k+1)}}\\
\frac{-q}{\sqrt{2(k+1)}}&0&\sqrt{\frac{k+1}{2}}&0\\
0&\frac{q}{\sqrt{2(k+1)}}&0&\sqrt{\frac{k+1}{2}}
\end{array}\right)\ .
\]
This matrix is again of class $C^1$ on $C_a$. In particular, it is non-singular and equal to the identity on $L$.
In fact, the matrix $\tilde{Q}$ can be written as $\tilde{Q}=\mathbb{1}+r_o^2\cdot\hat{Q}$, where $\hat{Q}$
is some matrix valued function on $C_a$ whose components are sums of functions of the form $f_l$ with $s_l\geq 0$. 
The spinor $\psi^o_{bc}$ is expressed with respect to the corresponding spinor frame $\tilde{h}_s$ by
\[ 
\psi^o_{bc}=[\ \tilde{h}_s\ ,\ \tilde{G}^{-1}\cdot\tilde{Q}^{-1}(-x_0b,x_0c,rb,rc)^\top\ ]\ .
\]
It holds
\[\Phi:=
\tilde{G}^{-1}\cdot\tilde{Q}^{-1}\left(\begin{array}{c}\!\!\!\!-x_0b\\x_0c\\rb\\rc\end{array}\right)
=
\tilde{G}^{-1}
\left(\begin{array}{ccc}
-\sqrt{\frac{k+1}{2}}\cdot x_0b&-&\frac{qrb}{\sqrt{2(k+1)}}\\[2mm]
\sqrt{\frac{k+1}{2}}\cdot x_0c&+&\frac{qrc}{\sqrt{2(k+1)}}\\[2mm]
\frac{qx_0b}{\sqrt{2(k+1)}}&+&\sqrt{\frac{k+1}{2}}\cdot rb\\[2mm]
\frac{qx_0c}{\sqrt{2(k+1)}}&+&\sqrt{\frac{k+1}{2}}\cdot rc
\end{array}\right)\ .
\]
From Lemma \ref{LM1} we know that the function $\frac{qx_0}{r}$ is of class $C^1$, since it behaves
like $r_o^2\cdot \frac{x_0}{r}$, whereby $\frac{x_0}{r}$ has order zero.
This observation is sufficient to conclude that the vector valued function $\Phi$ 
extends to a $C^1$-function  on $C_a$. Obviously, the extended $C^1$-function $\Phi$ is zero at the origin.
We can conclude that $\psi_{bc}^o$ extends to a $C^1$-spinor $\psi_{bc}$ on $\tilde{B}_a$ with zero at the origin.
If $(b,c)\neq 0$ the origin is the only zero of $\psi_{bc}$.
Moreover, since $\psi_{bc}^o$ is a twistor and $\psi_{bc}$ is $C^1$,
the expression
$\nabla^S_X\psi_{bc}+\frac{1}{n}X\cdot D^S\psi_{bc}$ is continuous on $\tilde{B}_a$ and zero on $\tilde{B}_a\smallsetminus \{0\}$
for all differentiable vector fields $X$. This shows that $\psi_{bc}$ satisfies the twistor equation in the origin.
Altogether we have proven yet the first three statements of Theorem \ref{TH1}.

The length square $u_{bc}$ of $\psi^>_{bc}=[e_s,(-x_0b,x_0c,rb,rc)^\top]$ is by definition (cf. \cite{Baum81}) equal to
\[
(\ \gamma_0\cdot (-x_0b,x_0c,rb,rc)^\top\ ,\ (-x_0b,x_0c,rb,rc)^\top\ )_{\CC^4}\ =\ (r^2-x_0^2)\cdot(b^2+c^2)\ .
\]
Obviously, the function $u_{bc}$ is smooth on $\tilde{B}_a$ and its zero set is $L_o$. We know already from 
Proposition \ref{PR2} that
$u_{bc}^{-2}\cdot g_a$ is a Ricci-flat metric on $\tilde{B}_a\smallsetminus L_o$, i.e., the function $u_{bc}$
provides a rescaling to an Einstein metric in the conformal class. It is well known that such a rescaling function
satisfies the partial differential equation 
\[
-u_{bc}\cdot Ric^0=(n-2)\cdot Hess(u_{bc})^0\ .
\]
It is interesting to note that the function $u_{bc}$ has a non-trivial zero set (cf. \cite{Gov04}). 

Furthermore, using the definition $g_a(V_{\psi_{bc}},X)=\langle\psi_{bc},X\cdot\psi_{bc}\rangle_S$ and calculating 
the products $(\gamma_0\cdot(-x_0b,x_0c,rb,rc)^\top,\gamma_i\cdot(-x_0b,x_0c,rb,rc)^\top)$ 
for $i=0,\ldots,4$ shows easily that
the spinor square of the twistor $\psi_{bc}$ is equal to
\[V_{\psi_{bc}}=(b^2+c^2)\cdot \big(-(x_0^2+r^2)e_0-2x_0re_1\big)=(b^2+c^2)\cdot V\ .\]
For $(b,c)\neq 0$ the vector field $V_{\psi_{bc}}$ is smooth with unique zero at the origin.
Finally, since $\alpha(V)=0$, we obtain $g_a(V,V)=-(r^2-x_0^2)^2$. This shows that $V_{\psi_{bc}}$ for $(b,c)\neq 0$ is
everywhere timelike except on $L_o$ where the spinor square is lightlike resp. zero only at the origin. \hfill$\Box$\\   

Corollary \ref{CO1} is a simple conclusion using Theorem \ref{TH1} and Proposition \ref{PR2}.

We add some remarks about the vector field $V$.
Although the twistor spinor $\psi_{bc}$ is not smooth, the spinor square $V_{\psi_{bc}}$ is smooth.
Since $\psi_{bc}$ is a twistor we immediately know that $V$ is a conformal Killing vector field for $g_a$ on $\tilde{B}_a$.
However, we simply reprove this statement here directly. Namely, it holds 
\[\begin{array}{ccl}
L_Vg_0&=& -4x_0\cdot g_0\\[2mm]
L_Vr^m&=&-2mx_0\cdot r^m\\[2mm]
L_Vr^2_o&=& \quad 0\\[2mm]
L_V\sigma_3^2&=& \quad 0\\[2mm]
L_V(-r^2(ar_o)^4\cdot \sigma_3^2)&=& -4x_0\cdot (-r^2(ar_o)^4\cdot \sigma_3^2)\\[2mm]
L_V\alpha&=&-4x_0\cdot \alpha\\[2mm]
L_V(a^4(r\beta)^{-2}r_o^2\cdot \alpha^2)&=& \ \  (4x_0-2\cdot 4\cdot x_0)(a^4(r\beta)^{-2}r_o^2\cdot \alpha^2)\\[1mm]
&=& -4x_0\cdot(a^4(r\beta)^{-2}r_o^2\cdot \alpha^2)\ .
\end{array}
\]
This proves that $L_Vg_a=-4x_0\cdot g_a$ on $\tilde{B}_a$, i.e., $V$ is a conformal Killing vector with 
$div^{g_a}(V)=-10\cdot x_0$.

An interesting property of $V$ is the fact that it is an essential conformal Killing vector field. 
The reason is that the spinor square $V_{D\psi_{bc}}$ of $D^S\psi_{bc}$ is given in $\{0\}\in\tilde{B}_a$ by 
\[
V_{D\psi_{bc}}\ =\ \frac{n}{2}\cdot grad(div(V_{\psi_{bc}}))\ ,
\]
which does not vanish, 
since it holds $D^S\psi_{bc}\neq 0$ in the origin.  
This argument is true for any metric in the conformal class $c_a=[g_a]$, 
i.e., the divergence of $V_{\psi_{bc}}$ resp. $V$ does not vanish identically with respect to any metric in $c_a$.

Finally, we want to state a reason why an extension of the metric $g_a$ on $B_a$ to $L$ with differentiable  
Weyl 
tensor 
has to be conformally flat
in order to preserve twistor spinors and the conformal Killing vector $V$. 
One observes the following facts. All integral curves of $V$ on $L$ converge in one flow direction to the origin, 
i.e., the origin is
in the closure of any integral curve on $L$. 
The length square $|W^{2,2}|^2$ of the
Weyl $(2,2)$-tensor is constant
along integral curves of $V$.
Moreover, with our assumptions we know that in the origin $W^{g_a}$ has to vanish (cf. \cite{Baum99}), i.e., $|W^{2,2}|^2$ is 
identically zero on
the closure $L$ of $L\smallsetminus L_o$.
Then, since $V$ inserted into the Weyl tensor $W^{g_a}$ produces zero
(cf. \cite{Baum99}) and $V$ is timelike on $L\smallsetminus L_o$, it follows that the length square of $W^{g_a}$ 
is non-negative on $L\smallsetminus L_o$ and it is zero if and only if the Weyl tensor vanishes. 
With the argument from before we can conclude that the Weyl tensor of the extension has to vanish on $L$.


\end{sloppypar}

\begin{thebibliography}{1111111}

\bibitem[Oba71]{Oba71} M. Obata. {\it The conjectures of conformal transformations of Riemannian manifolds.}  
Bull. Amer. Math. Soc.  77(1971), 265--270.

\bibitem[LF71]{LF71} J. Lelong-Ferrand. {\it Transformations conformes et quasi-conformes des varietes riemanniennes compactes.}
Acad. Roy. Belg. Cl. Sci. Mem. Coll. 8(2)  39,  no. 5, 44 pp. (1971).

\bibitem[Ale72]{Ale72}
D. Alekseevskii.
{\it Groups of conformal
transformations of Riemannian spaces}.
Mat. Sbornik 89(131) 1972 (in Russian),
English translation Math. USSR Sbornik
18(1972), 285-301.

\bibitem[Yos75]{Yos75}
Y. Yoshimatsu.
{\it On a theorem of Alekseevskii
concerning conformal transformations}. J. Math. Soc. Japan
28(1976), 278-289.

\bibitem[EH78]{EH78}
T. Eguchi, A.J. Hanson. {\it Asymptotically flat self-dual solutions to euclidean gravity.} Phys. Lett. B. 74(1978), 249-251.

\bibitem[Baum81]{Baum81} H. Baum.
{\it Spin-Strukturen und Dirac-Operatoren {\"u}ber pseudo-Riemannschen
Mannigfaltigkeiten.} No. 41 of {\it Teubner-Texte zur
Mathemtik}, Teubner, Leipzig
1981.


\bibitem[Lic88]{Lic88} A. Lichnerowicz. {\it Killing spinors, twistor
spinors and Hijazi inequality.} J. Geom. Phys. 5(1988), 2-18.

\bibitem[Lic89]{Lic89} A. Lichnerowicz. {\it On the twistor spinors.} Lett. Math. Phys. 18(1998), 333-345.

\bibitem[Lic90]{Lic90} A. Lichnerowicz. {\it Sur les zeros des spineurs-twisteurs.} C.R. Acad. Sci. Paris, Serie I, 310(1990), 19-22.

\bibitem[D'AG91]{DAG91} G. D'Ambra, M. Gromov. {\it Lectures on transformation groups: geometry and dynamics.}
Surveys in differential geometry (Cambridge, MA, 1990), 19-111, Lehigh Univ., Bethlehem, PA, 1991.

\bibitem[BFGK91]{BFGK91} H. Baum, Th. Friedrich,
R. Grunewald, I. Kath. {\it Twistor and Killing spinors on Riemannian
manifolds.} Teubner-Text No. 124, Teubner-Verlag, Stuttgart-Leipzig 1991.

\bibitem[KR95]{KR95}
W. K{\"u}hnel \& H.-B. Rademacher.
{\it Twistor spinors with
zeros}. Int. J. Math. 5(1994), 877-895.

\bibitem[KR96]{KR96}
W. K{\"u}hnel, H.-B. Rademacher.
{\it Twistor spinors and gravitational instantons.}
Lett. Math. Phys. 38(1996), 411-419.

\bibitem[KR98]{KR98} W. K{\"u}hnel, H.-B. Rademacher. {\it Asymptotically Euclidean manifolds and twistor spinors.}  Comm. Math. Phys.  196
(1998),  no. 1, 67--76.

\bibitem[Baum99]{Baum99} H. Baum. {\it Lorentzian twistor spinors and CR-geometry.}
J. Diff. Geom. and its Appl. 11(1999), no. 1, 69-96.



\bibitem[Lei99]{Lei99} F. Leitner. {\it Zeros of conformal vector fields and twistor spinors in Lorentzian geometry.}
SFB288 e-print no. 439, Berlin 1999.

\bibitem[Lei01]{Lei01} F. Leitner. {\it The twistor equation in Lorentzian geometry.} Dissertation HU Berlin, 2001.

\bibitem[Lei04]{Lei04} F. Leitner. {\it
A note on twistor spinors with zeros in Lorentzian geometry.}
e-print arXiv:math.DG/0406298, 2004.


\bibitem[Gov04]{Gov04} A.R. Gover. {\it Almost conformally Einstein manifolds and obstructions.}
electronic preprint, arXiv:math.DG/0412393, 2004.



\end{thebibliography}
\end{document}